\documentclass[11pt]{article}

\usepackage[top=2cm, bottom=3cm, left=2cm, right=2cm]{geometry}

\usepackage[abbrev]{amsrefs}

\makeatletter
\long\def\unmarkedfootnote#1{{\long\def\@makefntext##1{##1}\footnotetext{#1}}}
\makeatother

\usepackage[hidelinks,pdfencoding=auto]{hyperref}

\usepackage{graphicx}
\usepackage{amsmath, amssymb,
fancybox,enumerate,color,amsthm}
\usepackage{enumitem}

\begin{document}

\renewcommand{\theequation}{\arabic{section}.\arabic{equation}}

\newtheorem{theorem}{\bf Theorem}[section]
\newtheorem{lemma}{\bf Lemma}[section]
\newtheorem{proposition}{\bf Proposition}[section]
\newtheorem{corollary}{\bf Corollary}[section]
\newtheorem{remark}{\bf Remark}[section]
\newtheorem{example}{\bf Example}[section]
\newtheorem{definition}{\bf Definition}[section]

\newcommand{\supp}{\mathop{\mathrm{supp}}}
\newcommand{\PA}{\partial}
\newcommand{\ve}{\varepsilon}
\newcommand{\TR}{\textcolor{red}}
\newcommand{\TB}{\textcolor{blue}}
\newcommand{\ds}\displaystyle

\newcommand{\re}{{\Bbb R}}
\newcommand{\nid}{\noindent}

\newcommand{\nn}{\nonumber}
\newcommand{\hu}{\hat{u}}
\newcommand{\tw}{\widetilde{w}}
\newcommand{\D}{D_0^{1,2}}
\newcommand{\DK}{D_{0,K}^{1,2}}
\newcommand{\rhp}{\rightharpoonup}
\def\R{{\mathbb R}}
\newcommand{\la}{\langle}
\newcommand{\ra}{\rangle}
\newcommand{\ls}{\{}
\newcommand{\rs}{\}_{n=1,2\cdots}}
\newcommand{\hv}{\hat{v}}
\newcommand{\hvp}{\widehat{\vp}}
\newcommand{\hw}{\widehat{w}}
\newcommand{\cv}{\bar{v}}
\newenvironment{prf}[1]
   {{\noindent \bf Proof of {#1}.}}{\hfill \qed}
\newcommand{\eq}[1]{{\begin{equation}#1\end{equation}}}

\newcommand{\eps}{\varepsilon}
\newcommand{\dd}{\,d}
\newcommand{\abs}[1]{\lvert #1\rvert}
\newcommand{\norm}[1]{\lVert #1\rVert}
\newcommand{\old}[1]{\cite[#1]{FIRT25}}
\newcommand{\calF}{\mathcal F}
\newcommand{\calT}{\mathcal T}
\newcommand{\FIRT}{Fujishima--Ioku--Ruf--Terraneo}

\makeatletter
 \@addtoreset{equation}{section}
\makeatother
\def\theequation{\arabic{section}.\arabic{equation}}

\title{
A remark on singular solutions of semilinear elliptic equations in two dimensions
}
\author{\small
Yohei Fujishima\\
\small
Department of Mathematical and Systems Engineering,
Faculty of Engineering,
\\
\small
Shizuoka University, \\
\small
3-5-1 Johoku, Hamamatsu 432-8561, Japan, 
\vspace{5pt}\\
\quad\\
\small
Norisuke Ioku
\\
\small
Mathematical Institute, 
\\
\small
Tohoku University, \\
\small
Aramaki 6-3, Sendai 980-8578, Japan,
\vspace{5pt}\\
\quad\\
\small
Bernhard Ruf
\\
\small
Accademia di Scienze e Lettere - Istituto Lombardo,\\
\small
via Borgonuovo, 25, Milano 20121, Italy,
\\
\quad\\
{\small and}\\
\quad\\
\small
Elide Terraneo\footnote{Corresponding author: elide.terraneo@unimi.it (E. Terraneo)}
\\
\small
Dipartimento di Matematica ``F. Enriques'',\\
\small
Universit\`a degli Studi di Milano,\\
\small
via C. Saldini 50, Milano 20133, Italy.
}

\date{
}
\maketitle

\begin{abstract}
We study positive radial singular solutions of $-\Delta u=f(u)$ in two
space dimensions. Fujishima, Ioku, Ruf, and Terraneo constructed such
solutions under asymptotic assumptions on $f$ and an additional
convergence-rate condition, denoted by $(f_3)$. We prove their local
existence result, with the same approximate solution and remainder
estimate, under the original assumptions $(f_1)$ and $(f_2)$ alone. The proof
combines a contraction on finite intervals with a nonlinear energy
estimate uniform in the terminal point and requires only $C^2$
regularity of $f$ near infinity. We also obtain derivative asymptotics
and, under additional global assumptions, a singular Dirichlet solution
which solves the equation distributionally across the origin.
An analytic oscillatory example satisfies our hypotheses but violates
both $(f_3)$ and the derivative assumptions of Kikuchi and Kumagai.
For every fixed finite iterate of the exponential of order at least two,
we verify the hypotheses and quantify the remainder; condition $(f_3)$
fails for all iterates of order at least three.
\end{abstract}

\ \\
\noindent
{\small
\vspace{5pt}
\vspace{5pt}
}

\newpage
\section{Introduction}
\label{sec:main}
We consider positive radial singular solutions $u\in C^2(B_R\setminus \{0\})$ of
\begin{equation}\label{eq:1.1}
 -\Delta u=f(u)\quad\text{in }B_R\setminus\{0\}\subset\R^2,
 \qquad \lim_{|x|\to0}u(x)=\infty,
\end{equation}
where $B_R$ is the disk of radius $R$ centered at the origin. Our purpose
is to prove the existence of such solutions for a large class of exponential nonlinearities,
thereby
improving the previous results in 
 \cite{FIRT25}.

The construction of singular solutions in two dimensions is dictated
by the growth of the nonlinearity. For nonlinearities of exponential
type, explicit singular profiles provide useful starting points.
Ioku, Ruf, and Terraneo  \cite[Proposition~2.1]{IRT19} constructed a singular Dirichlet
solution for a nonlinearity equal to $s^{-3}e^{s^2}$ for large $s$; near
the origin its profile is exactly $\sqrt{-2\log |x|}$. Ibrahim, Kikuchi,
Nakanishi, and Wei \cite[Theorem~2.1]{IKNW21} constructed a singular Dirichlet solution
for $f(s)=s(e^{s^2}-1)$ and determined its asymptotic expansion. These
elliptic constructions were used in their studies of nonuniqueness for
heat equations.


Inspired by the explicit singular solution in 
\cite{IRT19}
one observes that there exist explicit singular solutions of $-\Delta v=g(v)$
for certain specific nonlinearities
$g$ (see \eqref{eq:1.5} and \eqref{eq:1.6}).
Using these nonlinearities and  singular solutions as a starting point in \cite{FIRT25}, we developed a construction of singular solutions for a class of general exponential nonlinearities $f$. Our method in \cite{FIRT25}  uses the
transformation introduced by Fujishima and Ioku \cite{FI18} to transfer
an explicit singular solution of a model equation with nonlinearity $g$ to an approximate
solution for $f$. The growth rate of the class of functions $f$ is described by a parameter $B$
defined below. A fixed-point argument then corrects the approximation,
provided that the rate of the discrepancy  of the nonlinearity $f$ from the model satisfies a smallness assumption (hypothesis $(f_3)$ in \cite{FIRT25}).
 This
additional hypothesis  excludes
some natural nonlinearities, including higher iterates of the
 exponential.
 In the present paper, using the same model approximation as in  \cite{FIRT25} we are able to extend our previous results to a larger class of functions, 
 without requiring that the rate of the discrepancy to the model satisfies $(f_3)$ in \cite{FIRT25}.

Kikuchi and Kumagai \cite[Theorem~1.1 and Appendix~A]{KK25}
subsequently constructed singular
solutions beyond the scope of $(f_3)$. They use a refined asymptotic
expansion under assumptions on $\gamma=\log f$, including
$\gamma\in C^5$ near infinity and bounds on the first three derivatives
of $\gamma\gamma''/(\gamma')^2$. Their examples include the triple and
quadruple exponential nonlinearities. 
In Proposition~\ref{prop:oscillatory} we give an analytic nonlinearity to which
our result applies, although both $(f_3)$ and the derivative assumptions
of \cite{KK25} fail.

We now state the assumptions and the main result. As in
\cite[Section~2]{FIRT25}, assume the following:
\begin{description}[style=nextline,leftmargin=2.9em,labelwidth=2.3em]
\item[$(f_1)$] For some $s_0>0$,
\[
 f\in C^2((s_0,\infty)),\quad f(s)>0,\quad f'(s)>0,\quad
 F(s):=\int_s^\infty\frac{\dd t}{f(t)}<\infty\quad(s>s_0).
\]
\item[$(f_2)$] Both limits exist and satisfy
\begin{equation*}
 \lim_{s\to\infty}f'(s)F(s)=1,\qquad
 B:=\lim_{s\to\infty}B_2[f](s)\in[1,\infty),
\end{equation*}
where
\begin{gather}
 \frac1{B_1[f](s)}
   =(-\log F(s))(1-f'(s)F(s)),\label{eq:1.2}\\
 \frac1{B_2[f](s)}
   =f'(s)F(s)(-\log F(s))^2
       \left(\frac{f(s)f''(s)F(s)}{f'(s)}-1\right).
       \label{eq:1.3}
\end{gather}
\end{description}
The parameter $B$ is finite, and the expressions involving $B_1$ and
$B_2$ are used for sufficiently large arguments.
 By l'H\^opital's rule and $(f_2)$,
\[
\lim_{s\to\infty}\frac{1-f'(s)F(s)}{1/(-\log F(s))}
=\lim_{s\to\infty}\frac1{B_2[f](s)}=\frac1B,
\]
 and so
	\begin{equation}\label{eq:1.4}B_1[f](s)\to B. 
	\end{equation}
	Moreover, as observed in \old{Remark~2.4},
	\begin{equation*}
		\qquad f'(s)F(s)<1
		\quad\text{for all sufficiently large $s$}.
		\end{equation*}
Let us now introduce the model nonlinearities of \old{(2.3)--(2.4)}.  We denote $r=|x|$, for $x\in \mathbb{R}^2$.

For $B>1$, 
put $B'=B/(B-1)$ and 
define
\begin{equation}\label{eq:1.5}
 g(s)=\frac4{BB'}s^{1-2B'}e^{s^{B'}},\qquad
 G(s)=\frac B4(s^{B'}+1)e^{-s^{B'}},\qquad 
 v(r)=\left(\log\frac1{r^2}\right)^{1/{B'}}.
\end{equation}
For $B=1$, define
\begin{equation}\label{eq:1.6}
 g(s)=4e^{e^s-2s},\qquad
 G(s)=\frac14(e^s+1)e^{-e^s},\qquad 
 v(r)=\log\log\frac1{r^2}.
\end{equation}
In both cases $G(s)=\int_s^\infty\dd t/g(t)$ and $-\Delta v=g(v)$
near the origin.
We define the approximate singular solution
\begin{equation}\label{eq:1.7}
\widetilde u(r)=F^{-1}(G(v(r))), \quad
 \text{where}\ G(v(r))=\frac B4r^2\left(\log\frac1{r^2}+1\right).
 \end{equation}
With $B_j[g]$ defined by \eqref{eq:1.2}--\eqref{eq:1.3} using $g,G$,
set
\begin{equation*}
 R_j(r)=\left|\frac1{B_j[f](\widetilde u(r))}
              -\frac1{B_j[g](v(r))}\right|,\qquad j=1,2.
\end{equation*}
Since $\lim_{s\to \infty}B_j[g](s)=B$ also holds for $g$ and $G$, \[ R_1(r)+R_2(r)\longrightarrow0
\quad(r\downarrow0).\] Remark that $\widetilde u(r)$ and $v(r)$ tend to infinity as $r\rightarrow 0$.
The additional hypothesis in \old{Theorem~2.1} is
\begin{description}[style=nextline,leftmargin=2.9em,labelwidth=2.3em]
\item[$(f_3)$] 
The discrepancies $R_1,R_2$ satisfy the following decay condition:
\begin{equation}
\label{eq:1.8}
\displaystyle  (-\log r)^{1/2}(R_1(r)+R_2(r))\longrightarrow0
 \quad(r\downarrow0).
 \end{equation}
\end{description}
Our main theorem removes this hypothesis.

\begin{theorem}\label{thm:main}
Assume $(f_1)$ and $(f_2)$. There are $r_0\in(0,1)$ and a positive radial
solution $u\in C^2(B_{r_0}\setminus\{0\})$ of \eqref{eq:1.1} such that
\begin{equation}\label{eq:1.9}
 u(r)=\widetilde u(r)+O\!\left(
 f(\widetilde u(r))F(\widetilde u(r))
 \sup_{0<t\le r}(R_1(t)+R_2(t))\right)
 \quad(r\downarrow0).
\end{equation}
Moreover, for all sufficiently small $r>0$,
\begin{equation}\label{eq:1.10}
 \left|\frac{u_r(r)}{\widetilde u_r(r)}-1\right|
 \le \frac C{\sqrt{1+\log(1/r^2)}}
       \sup_{0<t\le r}(R_1(t)+R_2(t)),
\end{equation}
and
\begin{equation}\label{eq:1.11}
 \frac{u_r(r)/f(u(r))}{v_r(r)/g(v(r))}\longrightarrow1
 \quad(r\downarrow0).
\end{equation}
\end{theorem}

The approximation and the remainder in \eqref{eq:1.9} agree with
\old{Theorem~2.1}. Thus $(f_1)$ and $(f_2)$ suffice for that conclusion, even
when the discrepancies $R_1,R_2$ tend to zero too slowly to satisfy
\eqref{eq:1.8} in $(f_3)$.
Under the additional assumptions on $[0,\infty)$ used
in \old{Theorem~2.2}, the local solution extends to a singular Dirichlet
solution.

\begin{corollary}\label{cor:global}
Assume $(f_1)$ and $(f_2)$, and suppose in addition that
\begin{equation}\label{eq:1.12}
 f\in C^1([0,\infty)),\quad
 f(s)>0,\quad f'(s)>0,\quad F(s)<\infty\quad(s>0).
\end{equation}
Then there exist $R>0$ and a positive radial function
$u_\infty\in C^2(B_R\setminus\{0\})
\cap C(\overline B_R\setminus\{0\})$ satisfying
\eqref{eq:1.9}--\eqref{eq:1.11} near the origin and
\[
 -\Delta u_\infty=f(u_\infty)\quad\text{in }B_R\setminus\{0\},
 \qquad u_\infty=0\quad\text{on }\partial B_R.
\]
Furthermore, $u_\infty,f(u_\infty)\in L^1(B_R)$ and
$-\Delta u_\infty=f(u_\infty)$ in $B_R$ in the sense of distributions.
Assumption $(f_3)$ is not required.
\end{corollary}

To make the distinction from \cite{KK25} explicit, write
\[
 \gamma(s)=\log f(s),\qquad
 H(s)=\frac{\gamma(s)\gamma''(s)}{(\gamma'(s))^2}.
\]
Condition (G2) in \cite[Section~1.2]{KK25} requires
$\gamma\in C^5$ near infinity and, for some $C>0$,
\begin{equation}\label{eq:1.13}
 |H^{(k)}(s)|\le C\left(\frac{\gamma'(s)}{\gamma(s)}\right)^k,
 \qquad k=1,2,3,
\end{equation}
for all sufficiently large $s$. The following example violates these
bounds despite being real analytic.

\begin{proposition}[An oscillatory analytic nonlinearity]
\label{prop:oscillatory}
The function
\begin{equation*}
 f(s)=\exp\!\left((s+2)^2+
               \frac{\sin((s+2)^4)}{(s+2)^9}\right),\qquad s\ge0,
\end{equation*}
is positive and strictly increasing on $[0,\infty)$, and extends to a
real analytic function on $(-2,\infty)$.
It satisfies $(f_1)$, $(f_2)$ with $B=2$, and
\eqref{eq:1.12}. However, it satisfies neither $(f_3)$ nor
the $k=1$ bound in \eqref{eq:1.13}.
Thus Theorem~\ref{thm:main} and Corollary~\ref{cor:global} apply to this
nonlinearity, which lies outside the hypotheses of both
\old{Theorem~2.1} and \cite[Theorem~1.1]{KK25}.
\end{proposition}



The strategy of the proof of  Theorem~\ref{thm:main} is the following.
		 Given $f$ satisfying $(f_1)$ and $(f_2)$,
consider the corresponding model function $g$ as  in \eqref{eq:1.5} or
 in \eqref{eq:1.6} with $B$ defined  in $(f_2)$.
Using the transformation introduced in \cite{FI18}, 
we have that
if $v$ is an explicit
 solution to $-\Delta v=g(v)$ near the origin then
$\tilde u:=F^{-1}(G(v))$ satisfies
\begin{equation}\label{eq:1.14}
-\Delta \tilde u=f(\tilde u)+\frac{|\nabla \tilde u|^2}{f(\tilde u)F(\tilde u)}\left(g'(v)G(v)-f'(\tilde u)F(\tilde u)\right).
\end{equation}
The idea is then to look for  a solution  to \eqref{eq:1.1} of the form
\begin{equation*}
	u=\tilde u +\theta
	\end{equation*}
where  $\theta $ is a correction term that, thanks to \eqref{eq:1.14}, satisfies
\begin{equation}\label{eq:1.15}
-\Delta  \theta=f(\tilde u+\theta)-f(\tilde u)-\frac{|\nabla \tilde u|^2}{f(\tilde u)F(\tilde u)}\left(g'(v)G(v)-f'(\tilde u)F(\tilde u)\right).
\end{equation}
In particular, we look for a correction term such that $|\theta|\leq f(\tilde u)F(\tilde u)$, 
as shown in
\cite{FIRT25}, because under this condition we can control uniformly the nonlinear term  $f(\tilde u+\theta)-f(\tilde u)$.
All the functions $\tilde u$, $v$ and $\theta$ are radially symmetric and they are defined with respect to the radial variable on an interval $(0,r_0)$ for some positive $r_0$.
Using the Emden-Fowler variable $\rho=1+\log\frac{1}{r^2}$ we consider the problem on $(\rho_0, \infty)$.
Then as in \cite{FIRT25}, by variation of parameters  we rewrite the differential equation in integral form on any interval $[\rho_0,M]$, with $M>\rho_0$, with  zero Cauchy data at $M$.
Now a pointwise truncation makes the nonlinear fixed point map well-defined and Lipschitz on each finite interval $[\rho_0,M]$ and an exponentially weighted norm gives a contraction (by the Bielecki renorming method for integral equations as in \cite{BP21}).
The fundamental ingredient is then an energy inequality satisfied by any  correction function $\theta$ such that  $|\theta|\leq f(\tilde u)F(\tilde u)$. This inequality makes the truncation inactive and  supplies the uniform control and decay needed to pass to  a locally convergent subsequence under $(f_1)$ and $(f_2)$. This limiting construction yields the solution on $[\rho_0,\infty)$.  Finally, since all the estimates in the proof contain the scalar factor $f(\tilde u)F(\tilde u)$,   we normalize the unknown correction $\theta$ dividing it by $f(\tilde u)F(\tilde u)$ to simplify computations.

As an application, we consider $f=\exp^{\circ n}$ for every fixed finite
integer $n\ge2$,
where 
$\exp^{\circ n}$ is the $n$-fold composition of the exponential function.
Singular solutions for iterated exponentials in
space dimensions at least three were studied by Ghergu and Goubet
\cite{GG20}. In two dimensions, the double exponential is covered by
\cite{FIRT25}, and the third and fourth iterates occur among the examples
of \cite{KK25}. Here we verify $(f_1)$ and $(f_2)$ for all finite $n\ge2$ and
compute the discrepancy from the  model. For $n\ge3$ it has order
$1/\log\rho$, which violates $(f_3)$, and Theorem~\ref{thm:main} yields an
explicit remainder in terms of iterated logarithms.

Section~\ref{sec:proof} contains the proofs of Theorem~\ref{thm:main}
and Corollary~\ref{cor:global}, together with the verification of
Proposition~\ref{prop:oscillatory}. Section~\ref{sec:examples} treats finite
iterates of the exponential.
Throughout this paper, we write $a \simeq b$ to mean that $a/b \to 1$
in the limit under consideration.

\section{Proofs of the main results}
\label{sec:proof}

We first collect estimates for the normalized correction equation and
prove an energy bound independent of the terminal point. We then
construct solutions on finite intervals by contraction and use this
bound to remove the truncation and pass to the infinite interval.

\subsection{The correction equation and coefficient estimates}
\label{subsec:normalized}

For $r>0$ sufficiently small, introduce the variables
\begin{equation*}
 \rho=1+\log\frac1{r^2},\quad
 \psi(\rho)=v(r),\quad \phi(\rho)=\widetilde u(r),\quad
 \eta(\rho)=\theta(r)=u(r)-\widetilde u(r),\quad \eps_j(\rho)=R_j(r).
\end{equation*}
Primes on functions of $\rho$ denote $\rho$ derivatives. Define
\begin{equation}\label{eq:2.1}
 h(s)=f(s)F(s),\qquad
 \delta(\rho)=\sup_{\tau\ge\rho}(\eps_1(\tau)+\eps_2(\tau)).
\end{equation}
Primes on $h$ denote derivatives with respect to $s$; primes on
$h(\phi)$ denote derivatives of the composite with respect to $\rho$.
Differentiating \eqref{eq:1.7} gives
\begin{equation}\label{eq:2.2}
 F(\phi)=G(\psi)=\frac B4\rho e^{1-\rho},\qquad
 \frac{\psi'}{g(\psi)G(\psi)}=1-\frac1\rho,\qquad
 \phi'=h(\phi)\left(1-\frac1\rho\right).
\end{equation}


The  transformation defined in \cite{FI18}, as applied in
\old{(3.2)--(3.5)}, gives the  equation \eqref{eq:1.15} for the correction term $\theta(r)$. In the variable  $\rho=1+\log\frac 1{r^2}$  with $\eta(\rho)=\theta(r)$ one may prove that \eqref{eq:1.15} reads as follows
\begin{equation}\label{eq:2.3}
 \eta''+\frac{e^{1-\rho}}4
           \bigl[f(\phi+\eta)-f(\phi)\bigr]+I=0,
\end{equation}
where the inhomogeneous term $I$ is defined  by
\begin{equation*}
	 I=h(\phi)\left(1-\frac1\rho\right)^2
	       \bigl[f'(\phi)F(\phi)-g'(\psi)G(\psi)\bigr].
	\end{equation*}
The profile $\phi$ satisfies the residual identity
$\phi''+e^{1-\rho}f(\phi)/4=I$. 
Hence a solution $\eta(\rho)=\theta(r)$ of
\eqref{eq:2.3} yields a solution of the original equation  $u(r)=\widetilde u(r)+\theta(r)$.
For a detailed derivation of equation \eqref{eq:2.3} for the correction term see \old{(3.2)--(3.5)}.

Set
\begin{align}
 L(\rho)&=-\frac3{16\rho^2}
           -\frac{1}{B\rho},\label{eq:2.4}\\
 N[\eta](\rho)&=\frac{e^{1-\rho}}4
       \bigl[f(\phi+\eta)-f(\phi)\bigr].\label{eq:2.5}
\end{align}
Then \eqref{eq:2.3} takes the form
\begin{equation}\label{eq:2.6}
 \eta''+\left(\frac1{B\rho}+\frac3{16\rho^2}\right)\eta
       +I+L\eta+N[\eta]=0.
\end{equation}

The functions $\rho^{1/4}\sin(2\sqrt{\rho/B})$ and
$\rho^{1/4}\cos(2\sqrt{\rho/B})$ solve the corresponding homogeneous
linear equation. As in \old{(3.14)}, for a finite terminal point
$M>\rho_0$, define
\begin{equation}\label{eq:2.7}
 T_M[\eta](\rho)=\sqrt B\int_\rho^M(\rho\tau)^{1/4}
   \sin\!\left(\frac2{\sqrt B}(\sqrt\rho-\sqrt\tau)\right)
   \bigl[I(\tau)+L(\tau)\eta(\tau)+N[\eta](\tau)\bigr]\dd\tau.
\end{equation}
A fixed point of $T_M$ satisfies \eqref{eq:2.6} and
$\eta(M)=\eta'(M)=0$.
We remark that the definitions of $L(\rho)$ and $ N[\eta](\rho)$ in \eqref{eq:2.4} and \eqref{eq:2.5} are different from the ones given in \cite{FIRT25}.
In particular the term  $L(\rho)$ vanishes to first order and not to second order as in \cite{FIRT25}. In fact,  since the integral operator $T_M$  is restricted  to a finite interval,  the second order vanishing property at infinity is not necessary.
By \old{Lemma~3.1}, whose proof uses only $(f_1)$ and $(f_2)$, the residual
satisfies
\begin{equation}\label{eq:2.8}
 |I(\rho)|\le C\frac{h(\phi)}\rho\eps_1(\rho),\qquad
 |I'(\rho)|\le C\frac{h(\phi)}{\rho^2}
                           (\eps_1(\rho)+\eps_2(\rho)).
\end{equation}
To estimate the nonlinear terms, we normalize the additive correction
by $\eta=h(\phi)y$. Put
\begin{equation}\label{eq:2.9}
 \tilde N(\rho,y)=\frac{f\Big(\phi(\rho)+h(\phi(\rho))y\Big)}{f(\phi(\rho))}-1,
 \qquad \tilde I(\rho)=\frac{B\rho I(\rho)}{h(\phi(\rho))}.
\end{equation}
The next lemma gives the coefficient estimates used in the energy
inequality and the contraction argument.

\begin{lemma}\label{lem:coefficients}
There are $\rho_*>1$ and constants $c,C>0$ such that, for $\rho\ge\rho_*$,
\begin{gather}
 \frac{(h(\phi))'}{h(\phi)}\le0,\qquad
   \left|\frac{(h(\phi))'}{h(\phi)}\right|\le\frac C\rho,\qquad
   \left|\frac{(h(\phi))''}{h(\phi)}\right|\le\frac C{\rho^2},\label{eq:2.10}\\
 |\tilde I(\rho)|+\rho|\tilde I'(\rho)|\le C\delta(\rho).\label{eq:2.11}
\end{gather}
For $|y|\le1$ the argument $\phi+h(\phi)y$ belongs to $(s_0,\infty)$ and
\begin{equation}\label{eq:2.12}
 \tilde N(\rho,0)=0,\qquad
 c\le\partial_y\tilde N(\rho,y)\le C,\qquad
 |\partial_\rho\tilde N(\rho,y)|\le\frac C{\rho^2}|y|.
\end{equation}
The last derivative is taken with $y$ fixed.
\end{lemma}

\begin{proof}
Direct differentiation and \eqref{eq:1.3} give $h'(s)=f'(s)F(s)-1$ and 
\begin{equation}\label{eq:2.13}
 h(s)h''(s)=f'(s)F(s)
       \left(\frac{f(s)f''(s)F(s)}{f'(s)}-1\right)
       =\frac1{(-\log F(s))^2B_2[f](s)}.
\end{equation}
Write
\begin{equation}\label{eq:2.14}
 -\log F(\phi)=\rho-1-\log\rho-\log(B/4)\simeq\rho.
\end{equation}
The first estimate below is given in \old{proofs of Lemmas~3.1 and~3.2};
the second follows from \eqref{eq:2.13}, \eqref{eq:2.14},
and $(f_2)$:
\begin{equation}\label{eq:2.15}
 h'(\phi)<0\ \text{and}\ 
h'(\phi)=O(\rho^{-1}),\qquad
 h(\phi)h''(\phi)=O(\rho^{-2}).
\end{equation}
Since $\phi'=h(\phi)(1-1/\rho)$, we have exactly
\begin{align}
 \frac{(h(\phi))'}{h(\phi)}&=h'(\phi)\left(1-\frac1\rho\right),\label{eq:2.16}\\
 \frac{(h(\phi))''}{h(\phi)}&=
  \bigl[h(\phi)h''(\phi)+h'(\phi)^2\bigr]
                      \left(1-\frac1\rho\right)^2
           +\frac{h'(\phi)}{\rho^2}.\label{eq:2.17}
\end{align}
Thus \eqref{eq:2.17} completes \eqref{eq:2.10}. In particular,
$h(\phi)$ is bounded above 
for large $\rho$ since $h(\phi)$ is decreasing and positive.  Since $\phi\to\infty$, we can
enlarge $\rho_*$ so that $\phi-h(\phi)>s_0$  for any $\rho\geq \rho_*$.
The first bound in \eqref{eq:2.11} follows from \eqref{eq:2.8}. For the
second bound use the exact identity
\[
 \tilde I'=\frac{BI}{h(\phi)}+\frac{B\rho I'}{h(\phi)}
       -\tilde I\frac{(h(\phi))'}{h(\phi)}
\]
and \eqref{eq:2.8}, \eqref{eq:2.10}.
We next prove \eqref{eq:2.12}. 
All statements in the rest of the proof are uniform
for $|y|\le1$. First remark that $\tilde{N}(\rho,0)=0$. 
In order to prove the remaining inequalities in \eqref{eq:2.12} we need some preliminary estimates.
First we prove that for any $s$ such that $|s-\phi|<h(\phi)$  the quotient $\frac{h(s)}{h(\phi)}$ is  bounded above and below by positive constants. Indeed,
 since $h'(s)\to0$ and $\phi-h(\phi)\to\infty$, the mean value
theorem gives
\begin{equation}\label{eq:2.18}
 \tfrac12h(\phi)\le h(s)\le2h(\phi)
 \quad\text{when }|s-\phi|\le h(\phi),
\end{equation}
after increasing $\rho_*$, if necessary.
Also for any $s$ such that $|s-\phi|<h(\phi)$ the ratio $\frac{f(s)}{f(\phi)}$ is bounded above and below by positive constants:
 indeed
\[
\log\frac{f(s)}{f(\phi)}
=\int_\phi^s\frac{f'(\sigma)F(\sigma)}{h(\sigma)}\dd \sigma,
\]
and \eqref{eq:2.18} and $f'F\to1$ bound the absolute value of
this integral.
The last needed estimate is  
the Taylor's theorem
to the second order of $h(\phi+h(\phi)y)$ centered at $\phi$, with coefficients that are uniformly controlled with respect to $|y|\leq 1$. 
First, we remark that for
 $|s-\phi|\le h(\phi)$, we have
$(-\log F)'(\sigma)=1/h(\sigma)$ and thanks to \eqref{eq:2.18} we get
\[
 |{-\log F(s)}+\log F(\phi)|
 \le \int_{\min\{s,\phi\}}^{\max\{s,\phi\}}\frac{\dd \sigma}{h(\sigma)}
 \le\frac{2|s-\phi|}{h(\phi)}\le2.
\]
Together with \eqref{eq:2.14}, this gives $-\log F(s)\simeq \rho$ on the whole interval
$[\phi-h(\phi),\phi+h(\phi)]$. The identity \eqref{eq:2.13}, \eqref{eq:2.18}, and $(f_2)$ now imply
\begin{equation}\label{eq:2.19}
 |h''(s)|\le\frac C{h(\phi)\rho^2}\qquad(|s-\phi|\le h(\phi)).
\end{equation}
Taylor's theorem therefore yields
\begin{equation}\label{eq:2.20}
 h(\phi+h(\phi)y)=h(\phi)\left[1+h'(\phi)y
                      +O\!\left(\frac{y^2}{\rho^2}\right)\right].
\end{equation}
We are now in a position to prove the remaining inequalities in \eqref{eq:2.12}.
Let $s=\phi+h(\phi)y$.  Direct differentiation gives
\begin{equation}\label{eq:2.21}
 \partial_y\tilde{N}
    =\frac{h(\phi)f'(s)}{f(\phi)}
    =\frac{h(\phi)}{h(s)}\,f'(s)F(s)\,\frac{f(s)}{f(\phi)}.
\end{equation}
By \eqref{eq:2.18}, $h(\phi)/h(s)$ is bounded above and below.
Together with $f'(s)F(s)\to1$ and the  bounds  for
$f(s)/f(\phi)$, this proves the middle estimate in \eqref{eq:2.12}.
Using \eqref{eq:2.2} and \eqref{eq:2.16}, direct differentiation of
\eqref{eq:2.9} gives
\begin{equation*}
 \partial_\rho\tilde{N}=
 \left(1-\frac1\rho\right)
 \left[(1+h'(\phi)y)\partial_y\tilde{N}
                  -(1+h'(\phi))(1+\tilde{N})\right].
\end{equation*}
Differentiating this identity with respect to $y$ gives
\begin{equation}\label{eq:2.22}
 \partial_y\partial_\rho\tilde{N}=
 \left(1-\frac1\rho\right)\partial_y\tilde{N}
 \left[(1+h'(\phi)y)h(\phi)\frac{f''(s)}{f'(s)}-1\right].
\end{equation}
To bound its last factor, write
\begin{equation*}
	(1+h'(\phi)y)h(\phi)\frac{f''(s)}{f'(s)}=	(1+h'(\phi)y)\frac{h(\phi)}{h(s)}\frac{h(s)f''(s)}{f'(s)}.
	\end{equation*}
Remark that by
 \eqref{eq:2.18}, \eqref{eq:2.19},
and $f'(s)F(s)\to1$ we get
\begin{equation}\label{eq:2.23}
 h(s)\frac{f''(s)}{f'(s)}
    =1+\frac{h(s)h''(s)}{f'(s)F(s)}
    =1+O(\rho^{-2}).
\end{equation}
On the other hand, by \eqref{eq:2.15}, \eqref{eq:2.20}
we obtain
\[
(1+h'(\phi)y)\frac{h(\phi)}{h(s)}
=\frac{1+h'(\phi)y}
{1+h'(\phi)y+O(y^2/\rho^2)}
=1+O(\rho^{-2}).
\]
Combining this last estimate and
\eqref{eq:2.23}, we finally get
\[
 (1+h'(\phi)y)\frac{h(\phi)}{h(s)}\frac{h(s)f''(s)}{f'(s)}
 =1+O(\rho^{-2}).
\]
Thus \eqref{eq:2.21} and \eqref{eq:2.22} imply
$|\partial_y\partial_\rho\tilde N|\le C\rho^{-2}$. Since
$\partial_\rho\tilde N(\rho,0)=0$, integration in $y$ proves the last
inequality in \eqref{eq:2.12}.
\end{proof}

\subsection{A uniform energy estimate}
\label{subsec:energy}

Substituting $\eta=h(\phi)y$ into \eqref{eq:2.3}, using
$e^{1-\rho}f(\phi)/4=h(\phi)/(B\rho)$, gives exactly
\begin{equation}\label{eq:2.24}
 y''+2\frac{(h(\phi))'}{h(\phi)} y'
       +\frac{\tilde N(\rho,y)+\tilde I}{B\rho}
       +\frac{(h(\phi))''}{h(\phi)} y=0.
\end{equation}
For the terminal conditions $\eta(M)=\eta'(M)=0$, the corresponding
conditions are $y(M)=y'(M)=0$. We derive a uniform estimate for
solutions which remain in $|y|\le1$.

Define the potential
\begin{equation*}
 V(\rho,y)=\int_0^y\tilde N(\rho,s)\dd s.
\end{equation*}
By Lemma~\ref{lem:coefficients}, for $\rho\ge\rho_*$ and $|y|\le1$,
\begin{equation}\label{eq:2.25}
 c y^2\le V(\rho,y)\le C y^2,\qquad
 y\tilde N(\rho,y)\ge c y^2,\qquad
 |\partial_\rho V(\rho,y)|\le C\rho^{-2}y^2.
\end{equation}

\begin{lemma}[Uniform terminal estimate]\label{lem:apriori}
There are $\rho_{**}\ge\rho_*$ and $C>0$, independent of the interval
endpoints, with the following property. If $M>\sigma\ge\rho_{**}$ and
$y\in C^2([\sigma,M])$ solves \eqref{eq:2.24}, satisfies
$y(M)=y'(M)=0$, and obeys $|y|\le1$ on $[\sigma,M]$, then
\begin{equation}\label{eq:2.26}
 |y(\rho)|+\sqrt\rho\,|y'(\rho)|\le C\delta(\rho)
       \qquad(\sigma\le\rho\le M).
\end{equation}
\end{lemma}

\begin{proof} The key idea to prove \eqref{eq:2.26} is to define an energy that bounds the quantity $|y(\rho)|+\sqrt\rho\,|y'(\rho)|$ and  whose evolution gives the desired estimates.
Choose a sufficiently large fixed $K>0$ and define, along a solution,
\begin{equation}\label{eq:2.27}
	E=\frac{B\rho}{2}(y')^2+V(\rho,y)+{\tilde I}y-\frac B4yy'+K{\tilde I}^2.
\end{equation}
See Remark \ref{R1} for a justification for this definition of  energy.
We first verify that this energy controls the needed quantities.
By \eqref{eq:2.25} and Young's inequality, $|{\tilde I}y|$ is bounded by
a small fixed multiple of $y^2$ plus $C {\tilde I}^2$. Also,
\[
\frac B4|yy'|\le\frac{B\rho}{4}(y')^2+\frac B{16\rho}y^2.
\]
Taking $K$ large and then $\rho_{**}$ large gives $c_E,C_E>0$ such that
\begin{equation}\label{eq:2.28}
	c_E\bigl[\rho(y')^2+y^2+{\tilde I}^2\bigr]
	\le E\le C_E\bigl[\rho(y')^2+y^2+{\tilde I}^2\bigr].
\end{equation}
Since $\partial_y V=\tilde N$, differentiating \eqref{eq:2.27} first gives
\[
	E'=\frac B4(y')^2
	+\left(B\rho y'-\frac B4y\right)y''
	+\partial_\rho V(\rho,y)
	+(\tilde N(\rho,y)+{\tilde I})y'+{\tilde I}'y+2K{\tilde I}{\tilde I}'.
\]
Substituting \eqref{eq:2.24} then gives
\begin{equation}\label{eq:2.29}
	\begin{split}
		E'&=\left(\frac B4-2B\rho\frac{(h(\phi))'}{h(\phi)}\right)(y')^2
		+\frac1{4\rho}y\tilde N(\rho,y)+\partial_\rho V(\rho,y)\\
		&\qquad +\left({\tilde I}'+\frac {\tilde I}{4\rho}\right)y+2K{\tilde I}{\tilde I}'
		-\left(B\rho y'-\frac B4y\right)\frac{(h(\phi))''}{h(\phi)} y
		+\frac B2\frac{(h(\phi))'}{h(\phi)} yy'.
	\end{split}
\end{equation}
By \eqref{eq:2.10} and \eqref{eq:2.25}, the first two terms of
\eqref{eq:2.29} are at least
$B(y')^2/4+c y^2/(4\rho)$.
For a fixed small $\epsilon>0$, \eqref{eq:2.11} implies
\begin{equation}\label{eq:2.30}
	\left|\left({\tilde I}'+\frac {\tilde I}{4\rho}\right)y\right|
	\le\frac\epsilon\rho y^2+\frac{C_\epsilon}\rho\delta^2,
	\qquad |2K{\tilde I}{\tilde I}'|\le\frac C\rho\delta^2.
\end{equation}
By \eqref{eq:2.10} and \eqref{eq:2.25}, the remaining errors satisfy
\begin{equation}\label{eq:2.31}
	\begin{split}
		|\partial_\rho V|
		+\left|\left(B\rho y'-\frac B4y\right)
		\frac{(h(\phi))''}{h(\phi)} y\right|
		+\left|\frac B2\frac{(h(\phi))'}{h(\phi)} yy'\right|
		&\le C\left(\frac{y^2}{\rho^2}
		+\frac{|yy'|}{\rho}\right) \\ 
		&\le\frac C{\rho^{3/2}}\bigl[\rho(y')^2+y^2\bigr],
	\end{split}
\end{equation}
where $|yy'|\le[\rho(y')^2+y^2]/(2\sqrt\rho)$ was used.
Combining this positive lower bound with \eqref{eq:2.30} and
\eqref{eq:2.31} in \eqref{eq:2.29}, choose $\epsilon$
small and then increase $\rho_{**}$. The error terms
$\epsilon y^2/\rho$ and
$C\rho^{-3/2}[\rho(y')^2+y^2]$ are then absorbed by
$B(y')^2/4+c y^2/(4\rho)$. Hence
\[
E'\ge\frac c\rho\bigl[\rho(y')^2+y^2\bigr]-\frac C\rho\delta^2.
\]
The upper bound in \eqref{eq:2.28} and ${\tilde I}^2\le C\delta^2$
give some $\mu>0$ such that
\begin{equation}\label{eq:2.32}
	E'\ge\frac\mu\rho E-\frac C\rho\delta(\rho)^2.
\end{equation}
All constants and thresholds chosen so far depend only on the tail
estimates for $f$, $B$, and the fixed $K$, not on $\sigma$ or $M$.

Multiplying
\eqref{eq:2.32} by $\rho^{-\mu}$ and integrating from
$\rho$ to $M$ yields
\begin{equation*}
	E(\rho)\le\left(\frac\rho M\right)^\mu E(M)
	+C\rho^\mu\int_\rho^M\frac{\delta(s)^2}{s^{\mu+1}}\dd s.
\end{equation*}
The terminal energy is
\begin{equation*}
	E(M)=K{\tilde I}(M)^2\le C\delta(M)^2.
\end{equation*}
Because $\delta$ is nonincreasing,
\[
\rho^\mu\int_\rho^M\frac{\delta(s)^2}{s^{\mu+1}}\dd s
\le\frac{\delta(\rho)^2}\mu,\qquad
\left(\frac\rho M\right)^\mu\delta(M)^2\le\delta(\rho)^2.
\]
It follows that $E(\rho)\le C\delta(\rho)^2$. The lower bound in
\eqref{eq:2.28} proves \eqref{eq:2.26}.
\end{proof}

\begin{remark}\label{R1}
We remark that a more  natural definition for the energy would be
	\[
	E_{0}=\frac{B\rho}{2}(y')^2+V(\rho,y)+\tilde Iy-\left(2+\frac B2\right)yy'+K\tilde I^2,
	\]
	for $K>0$ large enough. Indeed,
	multiplying the equation \eqref{eq:2.24} by $B\rho y'$ we get
	\[
	\begin{split}
		&B\rho y''y'+2B\rho\frac{(h(\phi))'}{h(\phi)} (y')^2
		+\tilde N(\rho,y)y'+\tilde Iy'
		+B\rho\frac{(h(\phi))''}{h(\phi)} yy'\\&=
		\left(\frac{B\rho}{2}(y')^2+V(\rho,y)+\tilde Iy-\left(2+\frac B2\right)yy'\right)'+R(\rho,y,y')=0
	\end{split}	
	\]
	where the remainder term $R(\rho,y,y')$ is negligible with respect to  $ y^2(\rho)+\rho (y'(\rho))^2
+\delta(\rho)^2
	$
 for large values of $\rho$ and we recall that 
	$2B\rho\frac{(h(\phi))'}{h(\phi)}
	\simeq
	 -2$.
	Moreover, by similar computations as in the case of the energy $E$ one can prove that  $E_0$
also
 controls 
 the needed quantities. Indeed,
	taking $K$ large and then $\rho_{**}$ large gives $c_{E_0},C_{E_0}>0$ such that
	\begin{equation*}
		c_{E_0}\bigl[\rho(y')^2+y^2+{\tilde I}^2\bigr]
		\le 
		E_0
		\le C_{E_0}\bigl[\rho(y')^2+y^2+{\tilde I}^2\bigr].
	\end{equation*}
	However, taking the derivative of $E_0$ and substituting $y''$ as in \eqref{eq:2.24}  we have
	\begin{equation*}
		\begin{split}
			E'_0&=\left(-2-2B\rho \frac{(h(\phi))'}{h(\phi)}\right)(y')^2+
			\left(\frac 2B+\frac 12\right)\frac y
			\rho\tilde {N}(\rho,y)+\partial_\rho V(\rho,y)\\
			&\qquad +\left(\tilde I'+\frac {(2+\frac B2)\tilde I}{B\rho}\right)y
			-\left(B\rho y'-(2+\frac B2)y\right)\frac{(h(\phi))''}{h(\phi)} y
			+(4+B)\frac{(h(\phi))'}{h(\phi)} yy' 
			+2K\tilde I \tilde I'.
		\end{split}
	\end{equation*}
	Since $\frac{(h(\phi))'}{h(\phi)}\simeq -\frac 1{B\rho}$, in this case it is not possible to bound from below, uniformly for large $\rho$, the derivative of $E_0$ by using the energy $E_0$ itself.
	A way to solve this difficulty is to replace in $E_0$, the term $(2+\frac B2) yy'$ by $(2+\frac B2-\epsilon)yy'$ for some $ \epsilon\in (0,2+\frac B2)$. In the definition of $E$  we replace $(2+\frac B2) yy'$  by $\frac B4yy'$, since this choice  simplifies 
	the computations. 
	\end{remark}
\subsection{Construction of solutions and completion of the proofs}
\label{subsec:contraction}

We now prove Theorem~\ref{thm:main} by constructing solutions on
finite intervals. Introduce the truncation
\begin{equation*}
 \pi(s)=\max\{-1,\min\{s,1\}\},\qquad s\in\R.
\end{equation*}
It satisfies $|\pi(s)|\le1$, $|\pi(s)-\pi(t)|\le|s-t|$, and
$\pi(s)=s$ for $|s|\le1$.

Fix $M>\rho_0\ge\rho_{**}$, where $\rho_0$ will be increased below
independently of $M$. On
\begin{equation}\label{eq:2.33}
 X_M=C([\rho_0,M];\R),\qquad
 \norm{y}_M=\sup_{\rho_0\le\rho\le M}e^{\rho-M}|y(\rho)|,
\end{equation}
define
\begin{equation}\label{eq:2.34}
 \widehat{\calT}_M[y](\rho)
       =\frac1{h(\phi(\rho))}T_M[h(\phi)\pi(y)](\rho),
 \qquad (h(\phi)\pi(y))(\tau)=h(\phi(\tau))\pi(y(\tau)).
\end{equation}
For each fixed finite $M$,
\begin{equation*}
 e^{\rho_0-M}\norm{y}_\infty\le\norm{y}_M\le\norm{y}_\infty,
\end{equation*}
so $X_M$ is complete in this norm. Since $h(\phi)>0$, the kernel is continuous,
and $\phi+h(\phi)\pi(y)>s_0$, the operator in
\eqref{eq:2.34} maps all of $X_M$ into $X_M$. In
particular, the nonlinear terms are well-defined on this space.

The weight in \eqref{eq:2.33} is the backward-time form of a Bielecki
weight (see \cite{BP21}). For the present map $\widehat{\calT}_M$, the
contraction estimate follows as below.

From \eqref{eq:2.2}, \eqref{eq:2.4}, and \eqref{eq:2.5}, for $|y|\le1$
we have the exact identity
\begin{equation*}
 Lh(\phi)y+N[h(\phi)y]
       =\frac{h(\phi)}{B\rho}\bigl[\tilde N(\rho,y)-y\bigr]
                          -\frac{3h(\phi)}{16\rho^2}y.
\end{equation*}
Thus Lemma~\ref{lem:coefficients} and the mean value theorem give
\begin{equation}\label{eq:2.35}
 |Lh(\phi)(y_1-y_2)+N[h(\phi)y_1]-N[h(\phi)y_2]|
       \le C\frac{h(\phi)}\rho|y_1-y_2|
       \qquad(|y_1|,|y_2|\le1).
\end{equation}
The constant is uniform for $\rho\ge\rho_{**}$.
For arbitrary $y_1,y_2\in X_M$, apply
\eqref{eq:2.35} to $\pi(y_1),\pi(y_2)$ in
\eqref{eq:2.7}. The inhomogeneous term $I$ cancels in the difference. Using
$h(\phi(\tau))\le h(\phi(\rho))$ for $\tau\ge\rho$ and the bound on the sine gives
\begin{equation*}
 \begin{split}
 |\widehat{\calT}_M[y_1](\rho)-\widehat{\calT}_M[y_2](\rho)|
 &\le C\int_\rho^M
   (\rho\tau)^{1/4}\frac{h(\phi(\tau))}{h(\phi(\rho))}\frac1\tau
                |y_1(\tau)-y_2(\tau)|\dd\tau\\
 &\le\frac C{\sqrt{\rho_0}}
                       \int_\rho^M|y_1(\tau)-y_2(\tau)|\dd\tau.
 \end{split}
\end{equation*}
Indeed, $\rho^{1/4}\tau^{-3/4}\le\rho^{-1/2}\le\rho_0^{-1/2}$.
Multiplying by $e^{\rho-M}$ and using
$|y_1(\tau)-y_2(\tau)|\le e^{M-\tau}\norm{y_1-y_2}_M$, we obtain
\[
 e^{\rho-M}|\widehat{\calT}_M[y_1](\rho)
                 -\widehat{\calT}_M[y_2](\rho)|
 \le\frac C{\sqrt{\rho_0}}\norm{y_1-y_2}_M
                                      \int_\rho^M e^{\rho-\tau}\dd\tau.
\]
The last integral is $1-e^{\rho-M}\le1$. Therefore
\begin{equation*}
 \norm{\widehat{\calT}_M[y_1]-\widehat{\calT}_M[y_2]}_M
       \le\frac C{\sqrt{\rho_0}}\norm{y_1-y_2}_M.
\end{equation*}
Choose $\rho_0$ so large that $C/\sqrt{\rho_0}\le1/2$. This choice is
independent of $M$. Banach's contraction principle on the whole space
$X_M$ gives a unique fixed point
\begin{equation*}
 y_M=\widehat{\calT}_M[y_M].
\end{equation*}
We record the equation and regularity supplied by this fixed point.
Let $\eta_M=h(\phi)y_M$. The finite integral in \eqref{eq:2.7} has a continuous
integrand. Its kernel vanishes on the diagonal, and its $\rho$ derivative
on the diagonal is $1$. Differentiation of that integral gives
\begin{equation}\label{eq:2.36}
 \begin{cases}
  \displaystyle 
 \eta_M''+
 \left(\frac1{B\rho}+\frac3{16\rho^2}\right)\eta_M
   +I+Lh(\phi)\pi(y_M)+N[h(\phi)\pi(y_M)]=0,\\[7pt]
 \eta_M(M)=\eta_M'(M)=0.
 \end{cases}
\end{equation}
Thus $\eta_M,y_M\in C^2([\rho_0,M])$ and $y_M(M)=y_M'(M)=0$.
Wherever $|y_M|\le1$, \eqref{eq:2.36} is exactly
\eqref{eq:2.6}, and $y_M$ solves \eqref{eq:2.24}.

We next use the energy estimate to show that the truncation is inactive.
Let $C_A$ be the constant in Lemma~\ref{lem:apriori}. Enlarge $\rho_0$
once more so that, in addition to the contraction condition,
\begin{equation*}
 C_A\delta(\rho_0)<\frac12.
\end{equation*}
All these choices precede the construction for an arbitrary $M>\rho_0$;
in particular, $\rho_0$ is fixed for the entire family.
Since $y_M(M)=0$, we have $|y_M|<1$ near $M$. If there is a point of
$[\rho_0,M]$ at which $|y_M|=1$, take the largest such point $\sigma$.
It exists by compactness, and $\sigma<M$. On $[\sigma,M]$ we have
$|y_M|\le1$, so the truncation is inactive there and
Lemma~\ref{lem:apriori} applies. At $\sigma$ it implies
\[
 1=|y_M(\sigma)|\le C_A\delta(\sigma)
                 \le C_A\delta(\rho_0)<\frac12,
\]
a contradiction. Continuity also precludes a value with $|y_M|>1$
without crossing $1$. Consequently $|y_M|<1$ throughout $[\rho_0,M]$,
so
\begin{equation*}
 \eta_M=T_M[\eta_M].
\end{equation*}
Lemma~\ref{lem:apriori}, now on the whole interval, yields
\begin{equation}\label{eq:2.37}
 |y_M(\rho)|+\sqrt\rho\,|y_M'(\rho)|\le C_A\delta(\rho)
       \qquad(\rho_0\le\rho\le M).
\end{equation}
It remains to pass to the infinite interval and recover the asserted
asymptotics.  Take $M_k\to\infty$.
On each fixed compact interval in
$[\rho_0,\infty)$, \eqref{eq:2.37} bounds $y_{M_k},y_{M_k}'$ and
\eqref{eq:2.24} bounds $y_{M_k}''$. A diagonal
Arzel\`a--Ascoli argument gives a subsequence converging in
$C^1_{\rm loc}([\rho_0,\infty))$ to $y$. The right side of the equation
for $y_{M_k}''$ then converges locally uniformly, so $y\in C^2$ solves
\eqref{eq:2.24}. Moreover,
\begin{equation}\label{eq:2.38}
 |y(\rho)|+\sqrt\rho\,|y'(\rho)|\le C\delta(\rho).
\end{equation}
Set $\eta=h(\phi)y$ and $u(r)=\phi(\rho)+\eta(\rho)$. Equation
\eqref{eq:2.3} holds, hence $-\Delta u=f(u)$ in a sufficiently small
punctured disk. The choice of $\rho_0$ gives $u>s_0>0$. Also $h(\phi)$ is bounded,
$\delta\to0$, and $\phi\to\infty$, so $\eta\to0$ and $u\to\infty$.
Finally,
\[
 |\eta(\rho)|\le Ch(\phi(\rho))\delta(\rho),\qquad
 \delta(\rho)=\sup_{\tau\ge\rho}(\eps_1(\tau)+\eps_2(\tau))
   =\sup_{0<t\le r}(R_1(t)+R_2(t)).
\]
This is \eqref{eq:1.9}. Using \eqref{eq:2.2} and
\eqref{eq:2.10} gives
\begin{equation*}
 \frac{\eta'}{\phi'}
   =\frac{((h(\phi))'/h(\phi))y+y'}{1-1/\rho}
   =O\!\left(\frac{\delta(\rho)}{\sqrt\rho}\right),
\end{equation*}
which proves \eqref{eq:1.10}. Also, \eqref{eq:2.12} and
\eqref{eq:2.38} imply
\[
 \frac{f(\phi+\eta)}{f(\phi)}=1+\tilde N(\rho,y)=1+O(\delta(\rho)).
\]
Using $\phi'/f(\phi)=\psi'/g(\psi)$, we conclude that
\[
 \frac{u_r(r)/f(u(r))}{v_r(r)/g(v(r))}
 =\frac{\phi'+\eta'}{\phi'}\,
   \frac{f(\phi)}{f(\phi+\eta)}\longrightarrow1.
\]
This proves \eqref{eq:1.11} and completes the proof of
Theorem~\ref{thm:main}.

\begin{proof}[Proof of Corollary~\ref{cor:global}]
Let $u$ be the solution given by Theorem~\ref{thm:main}. Since
$\widetilde u_r<0$, \eqref{eq:1.10} gives $u_r(r_1)<0$ for
some small $r_1>0$. Set $c=-r_1u_r(r_1)>0$. On any positive
continuation to the right of $r_1$, the radial equation implies
\[
 (ru_r)'=-rf(u)<0,\qquad
 ru_r(r)\le-c,\qquad
 u(r)\le u(r_1)-c\log\frac r{r_1}.
\]
Such a continuation cannot stop at a finite radius with a positive
limiting value: $u$ remains in $[0,u(r_1)]$, and the identity
\[
 ru_r(r)=r_1u_r(r_1)-\int_{r_1}^r t f(u(t))\dd t
\]
bounds $u_r$ on finite intervals with $r\ge r_1$. The ordinary
differential equation continuation theorem applies, and the logarithmic
upper bound forces a first zero at some
$R\le r_1\exp(u(r_1)/c)<\infty$. Joining this continuation to the local
solution gives the stated function $u_\infty$.

It remains to verify the equation in the nonpunctured ball. By
\eqref{eq:1.2}, \eqref{eq:1.4}, and \eqref{eq:2.16},
\[
 \frac{(h(\phi))'}{h(\phi)}=-\frac{1/B+o(1)}\rho.
\]
Integrating gives, for every $\epsilon>0$,
$h(\phi(\rho))\le C_\epsilon\rho^{-1/B+\epsilon}$ for large $\rho$;
compare \old{Lemma~3.5}. 
Using this property and  choosing
	$0<\epsilon<1/B$, one may proceed as in Lemma 3.6 in \cite{FIRT25} 
 to obtain
\begin{equation*}
 \begin{split}
 u_\infty(r)&=O((-\log r)^{1-1/B+\epsilon}),\qquad
 |r(u_\infty)_r(r)|=O((-\log r)^{-1/B+\epsilon}),\\
 f(u_\infty(r))&=O\!\left(\frac1{r^2(-\log r)^{1+1/B-\epsilon}}\right).
 \end{split}
\end{equation*}
These bounds imply
$u_\infty,f(u_\infty)\in L^1(B_R)$,
$r|(u_\infty)_r(r)|\to0$, and $r u_\infty(r)\to0$.
For a test function $\zeta\in C_c^\infty(B_R)$, integration by parts
on $B_R\setminus\overline B_r$ therefore gives
\[
 \int_{B_R}\bigl(u_\infty\Delta\zeta+f(u_\infty)\zeta\bigr)\dd x=0
\]
after letting $r\downarrow0$: the inner boundary terms are bounded by
$C_\zeta r\bigl(u_\infty(r)+|(u_\infty)_r(r)|\bigr)$ and tend to zero.
This proves the distributional assertion.
\end{proof}

\begin{proof}[Proof of Proposition~\ref{prop:oscillatory}]
Let $\gamma=\log f$. All asymptotic estimates in this proof are taken
as $s\to\infty$. Direct differentiation
gives
\begin{equation}\label{eq:2.39}
 \begin{split}
 \gamma'&=2(s+2)+4(s+2)^{-6}\cos((s+2)^4)
                    -9(s+2)^{-10}\sin((s+2)^4),\\
 \gamma''&=2-16(s+2)^{-3}\sin((s+2)^4)+O(s^{-7}),\\
 \gamma'''&=-64\cos((s+2)^4)+O(s^{-4}).
 \end{split}
\end{equation}
For $s\ge0$, we have
$\gamma'\ge2(s+2)-4(s+2)^{-6}-9(s+2)^{-10}>0$.
Thus $f,f'>0$ on $[0,\infty)$.
Real analyticity is immediate, and $\gamma(s)=(s+2)^2+O(s^{-9})$ implies
$F(s)<\infty$. Consequently $(f_1)$ and \eqref{eq:1.12}
hold.

To verify $(f_2)$, let $F_0(s)=\int_s^\infty e^{-(t+2)^2}\dd t$. The bound
$|(t+2)^{-9}\sin((t+2)^4)|\le(s+2)^{-9}$ for $t\ge s$ gives
$F(s)=F_0(s)(1+O(s^{-9}))$. Repeated integration by parts in the
Gaussian tail therefore yields
\begin{equation*}
 F(s)=\frac{e^{-(s+2)^2}}{2(s+2)}
       \left(1-\frac1{2(s+2)^2}+\frac3{4(s+2)^4}+O(s^{-6})\right).
\end{equation*}
We obtain
\begin{equation}\label{eq:2.40}
 \begin{split}
f'F &=1-\frac1{2(s+2)^2}+\frac3{4(s+2)^4}+O(s^{-6}),\\
-\log F &=(s+2)^2+\log(2(s+2))+O(s^{-2}).
 \end{split}
\end{equation}
Moreover, \eqref{eq:2.39} gives
\[
 \frac{\gamma''}{(\gamma')^2}
    =\frac1{2(s+2)^2}-\frac4{(s+2)^5}\sin((s+2)^4)+O(s^{-9}),
\]
and hence
\[
 \frac{ff''F}{f'}-1
  =f'F\left(1+\frac{\gamma''}{(\gamma')^2}\right)-1
  =\frac1{2(s+2)^4}-\frac4{(s+2)^5}\sin((s+2)^4)+O(s^{-6}).
\]
Substitution in \eqref{eq:1.3} proves
\begin{equation}\label{eq:2.41}
 \frac1{B_2[f](s)}
   =\frac12-\frac4{s+2}\sin((s+2)^4)
                   +O\!\left(\frac{\log s}{s^2}\right).
\end{equation}
Together with $f'F\to1$, this establishes $(f_2)$ with $B=2$.

We next check the failure of the derivative bound in \cite{KK25}.
Using the explicit derivatives in \eqref{eq:2.39}, rather
than differentiating an asymptotic remainder, we find
\begin{equation*}
 \begin{split}
 H&=\frac12-\frac4{(s+2)^3}\sin((s+2)^4)+O(s^{-7}),\\
 H'&=\frac{\gamma''}{\gamma'}
       +\frac{\gamma\gamma'''}{(\gamma')^2}
       -\frac{2\gamma(\gamma'')^2}{(\gamma')^3}\\
   &=-16\cos((s+2)^4)+O(s^{-4}),\qquad
 \frac{\gamma'}\gamma\simeq\frac2{s}.
 \end{split}
\end{equation*}
Along $s=(2\pi j)^{1/4}-2$, $j\to\infty$, we have $H'\to-16$, whereas
$\gamma'/\gamma\to0$. Thus \eqref{eq:1.13} fails already
for $k=1$. In fact, condition (G1) of \cite{KK25} does hold: $\gamma>0$
near infinity, $\int_s^\infty\dd t/\gamma(t)<\infty$, and $H\to1/2$.

Finally, the model for $B=2$ is
\[
 g(v)=v^{-3}e^{v^2},\qquad
 G(v)=\tfrac12(v^2+1)e^{-v^2},\qquad
 v(r)=\sqrt{-2\log r}.
\]
Direct substitution in \eqref{eq:1.3} gives
\[
 \frac1{B_2[g](v)}=\frac12
                    +O\!\left(\frac{\log v}{v^2}\right).
\]
Since $F(\widetilde u(r))=G(v(r))$,
\eqref{eq:2.40} implies $\widetilde u(r)+2\simeq v(r)$. Thus
\eqref{eq:2.41} and the definition of $R_2$ give
\begin{equation*}
 R_2(r)=\frac4{\widetilde u(r)+2}|\sin((\widetilde u(r)+2)^4)|
                   +O\!\left(\frac{\log\widetilde u(r)}{\widetilde u(r)^2}\right),
\end{equation*}
so
\[
 \sqrt{-\log r}\,R_2(r)
        =2\sqrt2\,|\sin((\widetilde u(r)+2)^4)|+o(1).
\]
As $r\downarrow0$, the continuous, strictly monotone function
$\widetilde u(r)+2$
runs through all sufficiently large values. The last expression
therefore does not tend to zero, which disproves $(f_3)$.
\end{proof}

\section{Finite iterates of the exponential}
\label{sec:examples}

Define
\begin{equation}\label{eq:3.1}
 E_0(s)=s,\qquad E_{k+1}(s)=\exp(E_k(s)).
\end{equation}
Thus $E_n=\exp^{\circ n}$ is the $n$-fold composition. The observation
that these nonlinearities have $B=1$ appears in \old{Section~2}. We
verify $(f_2)$ and compute the model discrepancies for every fixed finite
$n\ge2$.

\begin{proposition}[Iterated exponentials]\label{prop:iterated}
For every fixed integer $n\ge2$, the nonlinearity $f(s)=E_n(s)$ satisfies
$(f_1)$ and $(f_2)$ with $B=1$. Consequently, Theorem~\ref{thm:main} and
Corollary~\ref{cor:global} apply, with
\begin{equation}\label{eq:3.2}
 \widetilde u(r)=F^{-1}\!\left[
       \frac14r^2\left(\log\frac1{r^2}+1\right)\right].
\end{equation}
For $n\ge3$ the errors in \eqref{eq:2.1} satisfy
\begin{equation}\label{eq:3.3}
 \eps_j(\rho)=\frac{1+o(1)}{\log\rho}\quad(j=1,2),\qquad
 \delta(\rho)\simeq\frac2{\log\rho}.
\end{equation}
In particular, $(f_3)$ fails for every $n\ge3$. For these $n$, the solution
satisfies the explicit bound
\begin{equation}\label{eq:3.4}
 u(r)-\widetilde u(r)
  =O\!\left(\frac1{
     \rho(\log\rho)^2\displaystyle\prod_{j=2}^{n-2}\log^{\circ j}\rho}
       \right),\qquad \rho=1+\log\frac1{r^2},
\end{equation}
where an empty product is $1$. For every $n\ge2$,
\begin{equation}\label{eq:3.5}
 u(r)\simeq\log^{\circ(n-1)}\rho.
\end{equation}
\end{proposition}

\begin{proof}
Smoothness, positivity, and strict monotonicity are immediate from
\eqref{eq:3.1}. For $s\ge0$, $E_n(s)\ge e^s$, so $F(s)<\infty$.
In particular, both the tail assumptions $(f_1)$ and the additional global
hypotheses of Corollary~\ref{cor:global} hold.

To check $(f_2)$ with sufficient precision, put $t=E_{n-1}(s)$. The inverse
is $s=\log^{\circ(n-1)}t$ and its derivative is
\begin{equation*}
 Q_n(t)=\frac{\dd s}{\dd t}
       =\frac1{t\displaystyle\prod_{j=1}^{n-2}\log^{\circ j}t}.
\end{equation*}
All logarithms are taken at sufficiently large positive arguments.
Changing variables in the defining integral gives the exact formula
\begin{equation}\label{eq:3.6}
 F(s)=\int_t^\infty e^{-v}Q_n(v)\dd v.
\end{equation}
We apply the integration-by-parts argument of \old{Lemma~4.1} to this
representation.

Set
\begin{equation}\label{eq:3.7}
 d_n(t)=\sum_{j=1}^{n-2}
     \frac1{\displaystyle\prod_{k=1}^j\log^{\circ k}t}.
\end{equation}
Empty sums are zero, so $Q_2(t)=1/t$ and $d_2(t)=0$. Logarithmic
differentiation gives the exact identities
\begin{equation}\label{eq:3.8}
 {Q_n'}=-\frac{1+d_n}{t}{Q_n},\qquad
 {Q_n''}
    =\frac{(1+d_n)(2+d_n)-t d_n'}{t^2}{Q_n}.
\end{equation}
For $n\ge3$, the first summand in \eqref{eq:3.7} dominates, hence
$d_n(t)\simeq1/\log t$. Differentiating each product in \eqref{eq:3.7}
gives
\[
 t d_n'(t)=-\sum_{j=1}^{n-2}
 \frac1{\prod_{k=1}^j\log^{\circ k}t}
 \left(\sum_{k=1}^j\frac1{\prod_{\ell=1}^k\log^{\circ\ell}t}\right).
\]
Therefore $|t d_n'|\le d_n^2$. Differentiating this finite sum once more
shows $|t^2d_n''|\le C(d_n^2+d_n^3)$ on a sufficiently far tail.
The identities \eqref{eq:3.8} and differentiation of $Q_n^{''}$ once  more  then imply
\begin{equation}\label{eq:3.9}
 |Q_n^{(k)}(t)|\le C Q_n(t)t^{-k},\qquad k=1,2,3.
\end{equation}
The constants here and below may depend on $n$.

Integrating \eqref{eq:3.6} by parts three times yields
\begin{equation}\label{eq:3.10}
 F(s)=e^{-t}\bigl[Q_n(t)+Q_n'(t)+Q_n''(t)
                                  +O(Q_n(t)t^{-3})\bigr].
\end{equation}
Indeed, the remaining integral is
$\int_t^\infty e^{-v}Q_n'''(v)\dd v$. The function $Q_n$ is positive and
decreasing, so \eqref{eq:3.9} bounds its absolute value by
$C e^{-t}Q_n(t)t^{-3}$. All boundary terms at infinity vanish by
\eqref{eq:3.9} and the exponential factor.

In the next formulas all functions of $t$ are evaluated at
$t=E_{n-1}(s)$. Since $f(s)=e^t$ and $\dd t/\dd s=1/Q_n(t)$,
\begin{align}
 f'(s)F(s)
   &=1+\frac{Q_n'}{Q_n}+\frac{Q_n''}{Q_n}+O(t^{-3}),
       \label{eq:3.11}\\
 -\log F(s)&=t-\log Q_n(t)+O(t^{-1})=t+O(\log t),
       \label{eq:3.12}\\
 \frac{f(s)f''(s)F(s)}{f'(s)}-1
   &=\left(1-\frac{Q_n'}{Q_n}\right)f'(s)F(s)-1\nonumber\\
   &=\frac{Q_n''}{Q_n}-\left(\frac{Q_n'}{Q_n}\right)^2+O(t^{-3})
     =\frac{1+d_n-t d_n'}{t^2}+O(t^{-3}).
       \notag
\end{align}
Substituting these estimates into \eqref{eq:1.2} and \eqref{eq:1.3} gives
\begin{align}
 \frac1{B_1[f](s)}
   &=1+d_n(t)+O(\log t/t),\label{eq:3.13}\\
 \frac1{B_2[f](s)}
   &=1+d_n(t)-t d_n'(t)+O(\log t/t).\label{eq:3.14}
\end{align}
Equations \eqref{eq:3.11} and \eqref{eq:3.14} prove $(f_2)$ with
$B=1$, including $n=2$.

To compute the discrepancies, recall that for $B=1$ the model is
$g(s)=4e^{e^s-2s}$, with $\psi(\rho)=\log(\rho-1)$, and the computation in
\old{Example~4.4} gives
\begin{equation}\label{eq:3.15}
 \frac1{B_j[g](\psi(\rho))}
        =1+O(\log\rho/\rho),\qquad j=1,2.
\end{equation}
For $t=E_{n-1}(\phi(\rho))$, \eqref{eq:2.14} and
\eqref{eq:3.12} imply first $t\simeq\rho$ and then
\begin{equation}\label{eq:3.16}
 t=\rho+O(\log\rho).
\end{equation}
For $n\ge3$, we have $d_n(t)\simeq1/\log t$ and
$t d_n'(t)=O((\log t)^{-2})$. Subtracting \eqref{eq:3.15} from
\eqref{eq:3.13} and \eqref{eq:3.14} proves the two asymptotics for
$\eps_j$ in \eqref{eq:3.3}. Their tail supremum is asymptotic to
$2/\log\rho$: for every fixed small $\varepsilon>0$, the sum is between
$(2-\varepsilon)/\log\tau$ and $(2+\varepsilon)/\log\tau$ for every
sufficiently large $\tau$, and one takes the supremum over $\tau\ge\rho$.
It follows that
\[
 \sqrt\rho\,(\eps_1(\rho)+\eps_2(\rho))
       \simeq\frac{2\sqrt\rho}{\log\rho}\longrightarrow\infty.
\]
Since $-\log r=(\rho-1)/2$, this is exactly the failure of $(f_3)$.

Finally, \eqref{eq:3.10} and \eqref{eq:3.16} give
\[
 h(\phi(\rho))=f(\phi)F(\phi)\simeq Q_n(\rho)
   =\frac1{\rho\prod_{j=1}^{n-2}\log^{\circ j}\rho}.
\]
Theorem~\ref{thm:main}, together with \eqref{eq:3.3}, now gives
\eqref{eq:3.4}. Applying the increasing inverse
$\log^{\circ(n-1)}$ to \eqref{eq:3.16} gives
$\phi(\rho)\simeq\log^{\circ(n-1)}\rho$; the additive error tends to zero.
This proves \eqref{eq:3.5}. For $n=2$, the same leading
asymptotic follows from \eqref{eq:3.16} and Theorem~\ref{thm:main}.
\end{proof}

\begin{remark}
The double exponential $E_2(s)=e^{e^s}$ already satisfies $(f_3)$, as shown in
\old{Example~4.4}. In that case
$u-\widetilde u=O(\log\rho/\rho^2)$, consistently with that example.
For $n\ge3$, the estimate \eqref{eq:3.4} applies despite
the failure of $(f_3)$. For instance, $e^{e^{e^s}}$ and $e^{e^{e^{e^s}}}$ have the remainders
$O(1/[\rho(\log\rho)^2])$ and
$O(1/[\rho(\log\rho)^2\log\log\rho])$, respectively, relative to
\eqref{eq:3.2}. For $n=1$, $f(s)=e^s$ gives $F(s)=e^{-s}$ and
$1/B_1[f]=1/B_2[f]=0$, so the required finite $B$ in $(f_2)$ is absent.
Thus the case $n=1$ lies outside the hypotheses of Theorem~\ref{thm:main}.
\end{remark}

\noindent 
\textbf{Declaration of AI Use.} 
The authors consulted ChatGPT (OpenAI, GPT-6 Astra and GPT-6.1 Sol
models) to explore whether the assumption $(f_3)$ could be removed from the
main result in \cite{FIRT25}. In response, ChatGPT suggested an approach
to removing this assumption by focusing on the Bielecki renorming
technique discussed in \cite{BP21} and the energy $E$. The authors
subsequently developed and refined the suggested approach and
independently verified the resulting argument. The authors take full
responsibility for the correctness and presentation of the final results.

\noindent 
\textbf{Acknowledgements.} 
This work was partially funded by JSPS KAKENHI (grant numbers: 23K03179, 26K00610). 
The research was partially supported by the Italian  Ministry of University and Research (MUR) through the ``Dipartimenti di Eccellenza'' Programme 2023-2027 (CUP G43C22004580005) -- Dipartimento di Matematica dell'Universit\`a degli Studi di Milano
and through the MUR project PRIN2022\allowbreak PNRR No.~P2022Z7ZAJ.
This paper was partially written during a research visit of the last author to the Department of Mathematics at Tohoku University. The last author gratefully acknowledges the 
department's
warm hospitality.

\end{document}